\documentclass[11pt,doublespace]{article}
\begin{document}


\newtheorem{thm}{Theorem}
\newtheorem{lem}[thm]{Lemma}
\newtheorem{cor}[thm]{Corollary}
\newtheorem{ex}[thm]{Example}
\newtheorem{prop}[thm]{Proposition}
\newtheorem{remark}[thm]{Remark}
\newtheorem{coun}[thm]{Counterexample}
\newtheorem{defn}[thm]{Definition}
\newtheorem{conj}{Conjecture}
\newtheorem{problem}{Problem}
\newcommand\ack{\section*{Acknowledgement.}}

\newcommand{\etal}{{\it et al. }}

\newcommand{\bbP}{{\rm I\hspace{-0.8mm}P}}
\newcommand{\bbE}{{\rm I\hspace{-0.8mm}E}}
\newcommand{\bbF}{{\rm I\hspace{-0.8mm}F}}
\newcommand{\bbI}{{\rm I\hspace{-0.8mm}I}}
\newcommand{\bbR}{{\rm I\hspace{-0.8mm}R}}
\newcommand{\bbRp}{{\rm I\hspace{-0.8mm}R}_+}
\newcommand{\bbN}{{\rm I\hspace{-0.8mm}N}}
\newcommand{\bbC}{{\rm C\hspace{-2.2mm}|\hspace{1.2mm}}}
\newcommand{\bbD}{{\rm I\hspace{-0.8mm}D}}
\newcommand{\bbQ}{\bf Q}
\newcommand{\bbZ}{{\rm \rlap Z\kern 2.2pt Z}}
\newcommand{\bbK}{{\rm I\hspace{-0.8mm}K}}

\newcommand{\matP}{{\bbP}}
\newcommand{\mattildeP}{\tilde{\bbP}}
\newcommand{\matPN}[1]{{\bbP}_{#1}}
\newcommand{\matPP}[1]{{\bbP}_{#1}^0}
\newcommand{\matE}{{\bbE}}
\newcommand{\mattildeE}{\tilde{\bbE}}
\newcommand{\matEP}[1]{{\bbE}_{#1}^0}
\newcommand{\matF}{{\bbF}}
\newcommand{\matR}{{\bbR}}
\newcommand{\matRp}{{\bbRp}}
\newcommand{\matN}{{\bbN}}
\newcommand{\matZ}{{\bbZ}}
\newcommand{\matI}{{\bbI}}
\newcommand{\matK}{{\bbK}}
\newcommand{\matQ}{{\bbQ}}
\newcommand{\matC}{{\bbC}}
\newcommand{\matD}{{\bbD}}

\newcommand{\calL}{{\cal L}}
\newcommand{\calM}{{\cal M}}
\newcommand{\calN}{{\cal N}}
\newcommand{\calF}{{\cal F}}
\newcommand{\calG}{{\cal G}}
\newcommand{\calD}{{\cal D}}
\newcommand{\calB}{{\cal B}}
\newcommand{\calH}{{\cal H}}
\newcommand{\calI}{{\cal I}}
\newcommand{\calP}{{\cal P}}
\newcommand{\calQ}{{\cal Q}}
\newcommand{\calS}{{\cal S}}
\newcommand{\calT}{{\cal T}}
\newcommand{\calC}{{\cal C}}
\newcommand{\calK}{{\cal K}}
\newcommand{\calX}{{\cal X}}
\newcommand{\cals}{{\cal S}}
\newcommand{\calE}{{\cal E}}

\newcommand{\koniecmat}{\,}

\newcommand{\eqd}{\ =_{\rm d}\ }
\newcommand{\toto}{\leftrightarrow}
\newcommand{\eqdistr}{\stackrel{\rm d}{=}}
\newcommand{\as}{\stackrel{\rm a.s.}{=}}
\newcommand{\convdistr}{\stackrel{\rm d}{\rightarrow}}
\newcommand{\convweak}{{\Rightarrow}}
\newcommand{\convas}{\stackrel{\rm a.s.}{\rightarrow}}
\newcommand{\convfidi}{\stackrel{\rm fidi}{\rightarrow}}
\newcommand{\convprob}{\stackrel{p}{\rightarrow}}
\newcommand{\behaveslikedistr}{\stackrel{d}{\sim}}
\newcommand{\deff}{\stackrel{\rm def}{=}}
\newcommand{\bis}{{'}{'}}
\newcommand{\Cov}{{\rm Cov}}
\newcommand{\Var}{{\rm Var}}
\newcommand{\Exp}{{\rm E}}

\newcommand{\nd}{n^{\delta}}
\newcommand{\koniec}{\newline\vspace{3mm}\hfill $\odot$}

\newcommand{\longmbox}[1]{\qquad {\mbox{\rm {#1}}}\qquad }
\newcommand{\shortmbox}[1]{\; {\mbox{\rm {#1}}}\; }
\newcommand{\truncatedintegral}{\int_{1-\frac{k_n}{n}}^{1-\frac{1}{n}}}


\title{Sums of extreme values of subordinated long-range dependent sequences: moving averages with finite variance  \protect}
\author{Rafa{\l} Kulik\thanks{Research supported in part by NSERC Canada Discovery Grants of Mikl\'{o}s Cs\"{o}rg\H{o}, Donald Dawson and Barbara Szyszkowicz at Carleton University}\\
University of Sydney \\
School of Mathematics and Statistics \\
F07, University of Sydney\\
NSW 2006, Australia
\\ rkuli@maths.usyd.edu.au\\
and\\
Mathematical Institute\\
University of Wroc{\l}aw\\
Pl. Grunwaldzki 2/4\\
50-384 Wroc{\l}aw\\
Poland}

\maketitle

\begin{abstract}
In this paper we characterize the limiting behavior of sums of
extreme values of long range dependent sequences defined as
functionals of linear processes with finite variance. The extremal
sums behave completely different by compared to the i.i.d case. In
particular, though we still have asymptotic normality, the scaling
factor is relatively bigger than in the i.i.d case, meaning that the
maximal terms have relatively smaller contribution to the whole sum.
Also, the scaling need not depend on the tail index of the
underlying marginal distribution, as it is well-known to be so in
the i.i.d. situation. Furthermore, subordination may completely
change the asymptotic properties of sums of extremes.
\end{abstract}



\\
\noindent{\bf Running title:} Sums of extremes and LRD
\section{Introduction}
Let $\{\epsilon_i,i\ge 1\}$ be a centered sequence of i.i.d. random
variables. Consider the class of stationary linear processes
\begin{equation}\label{model}
X_i=\sum_{k=0}^{\infty}c_k\epsilon_{i-k},\ \ \ i\ge 1 .
\end{equation}
We assume that the sequence $c_k$, $k\ge 0$, is regularly varying
with index $-\beta$, $\beta\in (1/2,1)$. This means that $c_k\sim
k^{-\beta}L_0(k)$ as $k\to\infty$, where $L_0$ is slowly varying at
infinity. We shall refer to all such models as long range dependent
(LRD) linear processes. In particular, if the variance exists (which
is assumed throughout the whole paper), then the covariances
$\rho_k:=\Exp X_0X_k$ decay at the hyperbolic rate,
$\rho_k=k^{-(2\beta-1)}L(k)$, where
$\lim_{k\to\infty}L(k)/L_0^2(k)=B(2\beta-1,1-\beta)$ and
$B(\cdot,\cdot)$ is the beta-function. Consequently, the covariances
are not summable (cf. \cite{GiraitisSurgailis2002}).

Assume that $X_1$ has a continuous distribution function $F$. For
$y\in (0,1)$ define $Q(y)=\inf\{x:F(x)\ge y\}=\inf\{x:F(x)= y\}$,
the corresponding (continuous) quantile function. Given the ordered
sample $X_{1:n}\le\cdots\le X_{n:n}$ of $X_1,\ldots,X_n$, let
$F_n(x)=n^{-1}\sum_{i=1}^n1_{\{X_i\le x\}}$ be the empirical
distribution function and $Q_n(\cdot)$ be the corresponding
left-continuous sample quantile function, i.e. $Q_n(y)=X_{k:n}$ for
$\frac{k-1}{n}<y\le \frac{k}{n}$. Define $U_i=F(X_i)$ and
$E_n(x)=n^{-1}\sum_{i=1}^n1_{\{U_i\le x\}}$, the associated uniform
empirical distribution. Denote by $U_n(\cdot)$ the corresponding
uniform sample quantile function.\\

Assume that $E\epsilon_1^2<\infty$. Let $r$ be an integer and define
$$
Y_{n,r}=\sum_{i=1}^n\sum_{1\le j_1<\cdots\le
j_r}\prod_{s=1}^rc_{j_s}\epsilon_{i-j_s},\qquad n\ge 1,
$$
so that $Y_{n,0}=n$, and $Y_{n,1}=\sum_{i=1}^nX_i$. If
$p<(2\beta-1)^{-1}$, then
\begin{equation}\label{eq-varance-behaviour}
\sigma_{n,p}^2:=\Var (Y_{n,p})\sim n^{2-p(2\beta-1)}L_0^{2p}(n).
\end{equation}
Define now the general empirical, the uniform empirical, the general
quantile and the uniform quantile processes respectively as follows:
$$
\beta_n(x)=\sigma_{n,1}^{-1}n(F_{n}(x)-F(x)),\qquad x\in {\matR},
$$
$$
\alpha_n(y)=\sigma_{n,1}^{-1}n(E_{n}(y)-y),\qquad y\in (0,1),
$$
$$
q_n(y)=\sigma_{n,1}^{-1}n(Q(y)-Q_{n}(y)),\qquad y\in (0,1),
$$
$$
u_n(y)=\sigma_{n,1}^{-1}n(y-U_{n}(y)),\qquad y\in (0,1).
$$
The aim of this paper is to study the asymptotic behavior of trimmed
sums based on the ordered sample $X_{1:n}\le\cdots\le X_{n:n}$
coming from the long range dependent sequence defined by
(\ref{model}).

Let $T_n(m,k)=\sum_{i=m+1}^{n-k}X_{i:n}$ and note that (see below
for a convention concerning integrals)
\begin{equation}\label{eq-trimmedmeans-represenation}
T_n(m,k)=n\int_{m/n}^{1-k/n}Q_n(y)dy .
\end{equation}
Ho and Hsing observed in \cite{HoHsing} that, under appropriate
conditions on $F$, as $n\to\infty$,
\begin{equation}\label{eq-HoHsing-quantiles}
\sup_{y\in
[y_0,y_1]}\left|q_n(y)+\sigma_{n,1}^{-1}\sum_{i=1}^nX_i\right|=o_P(1),
\end{equation}
where $0<y_0<y_1<1$. Equation (\ref{eq-HoHsing-quantiles}) means
that, in principle, {\it the quantile process can be approximated by
partial sums, independently of $y$}. This observation, together with
(\ref{eq-trimmedmeans-represenation}), yields the asymptotic
normality of the trimmed sums in case of heavy trimming
$m=m_n=[\delta_1 n]$, $k=k_n=[\delta_2 n]$, where
$0<\delta_1<\delta_2<1$ and $[\cdot]$ is the integer part (see
\cite[Corollary 5.2]{HoHsing}). This agrees with the i.i.d.
situation (see \cite{Stigler}).

However, the representation (\ref{eq-trimmedmeans-represenation})
requires some additional assumptions on $F$. In order to avoid them,
we may study asymptotics for the trimmed sums via the integrals of
the form $\int\alpha_n(y)dQ(y)$. This approach was initiated in two
beautiful papers by M. Cs\"{o}rg\H{o}, S. Cs\"{o}rg\H{o},
Horv\'{a}th and Mason, \cite{CsorgoCsorgoHorvathMason-1986a},
\cite{CsorgoCsorgoHorvathMason-1986b}. Then, S. Cs\"{o}rg\H{o},
Haeusler, Horv\'{a}th and Mason took this route to provide the full
description of the weak asymptotic behavior of the trimmed sums in
the i.i.d. case. The list of the papers written by these authors on
this particular topic is just about as long as this introduction.
Therefore we refer to \cite{CsorgoHaeusslerMason} for an extensive
up-to-date discussion and a survey of results.

In the LRD case, instead of using the Brownian bridge approximation,
we can use the reduction principle for the general empirical
processes as studied in \cite{GiraitisSurgailis2002},
\cite{HoHsing}, \cite{Surgailis-Koul-2002} or \cite{Wu2003} (see
Lemma \ref{thm-HoHsing} below). We can then use an approach that is
similar to that the above mentioned authors to establish asymptotic
normality in case of light, moderate and heavy trimming with the
scaling factor $\sigma_{n,1}^{-1}$, which is the same as for the
whole partial sum. So, in this context the situation is similar to
the i.i.d. case and for details we refer the
reader to the technical report \cite{Kulik-OuldHaye}.\\

The most interesting phenomena, however, occur when one deals with
the $k_n$-extreme sums, $\sum_{i=n-k_n+1}^nX_i$. If $F(0)=0$ and
$1-F(x)=x^{-\alpha}$, $\alpha>2$, then in the i.i.d situation we
have $$a_n\sum_{i=n-k_n+1}^nX_i-c_n\convdistr Z,$$ where the scaling
factor is $a_n=\left(nk_n^{-1}\right)^{1/2-1/\alpha}n^{-1/2}$, $c_n$
is a centering sequence and $Z$ is a standard normal random variable
(see \cite{CsorgoMason1986}). In the LRD case we still obtain
asymptotic normality. However, although the Ho and Hsing result
(\ref{eq-HoHsing-quantiles}) does not say anything about the
behavior of the quantile process in the neighborhood of $0$ and $1$,
the somewhat imprecise statement that {\it the quantile process can
be approximated by partial sums, independently of $y$} suggests that
\begin{itemize}
\item
a required scaling factor would not depend on the tail index
$\alpha$.
\end{itemize}
Indeed, we will show in Theorem \ref{thm-extreme-sums} that the
appropriate scaling in case $1-F(x)=x^{-\alpha}$ is
$(nk_n^{-1})\sigma_{n,1}^{-1}$. Removing the scaling for the whole
sums ($n^{-1/2}$ and $\sigma_{n,1}^{-1}$ in the i.i.d. and LRD
cases, respectively), we also see that
\begin{itemize}
\item the scaling in the LRD situation is greater, meaning that
the $k_n$-extreme sums contribute relatively {\it less} to the whole
sum compared to the i.i.d situation. This also is quite intuitive.
Since the dependence is very strong, it is very unlikely that we
have few big observations, which is a typical case in the i.i.d.
situation. Rather, if we have one big value, we have a lot of them.
\end{itemize}
One may ask, whether such phenomena are typical for all LRD
sequences. Not likely. Define $Y_i=G(X_i)$, $i\ge 1$, with some
real-valued measurable function $G$. In particular, taking
$G=F_Y^{-1}F$ we may obtain a LRD sequence with the arbitrary
marginal distribution function $F_Y$. Assume for a while that $F$,
the distribution of $X_1$, is standard normal and that $q_n(\cdot)$
is the quantile process associated with the sequence $\{Y_i,i\ge
1\}$. Following \cite{CSW2006} we observed in \cite[Section
2.2]{CsorgoKulik2006a} and \cite{CsorgoKulik2006b} that $q_n(\cdot)$
is, up to a constant, approximated by
$\phi(\Phi^{-1}(y))/f_Y(F_Y^{-1}(y))\sigma_{n,1}^{-1}\sum_{i=1}^nX_i$.
Here, $f_Y$ is the density of $F_Y$ and $\phi$, $\Phi$ are the
standard normal density and distribution, respectively. In the
non-subordinated case, $Y_i=X_i$, and the factor
$\phi(\Phi^{-1}(y))/f_Y(F_Y^{-1}(y))$ disappears. Nevertheless, from
this discussion it should be clear that the limiting behavior of the
extreme sums in the subordinated case $Y_i=G(X_i)$ is different,
namely (see Theorem \ref{thm-extreme-sums})
\begin{itemize}
\item the scaling depends on the marginal distributions of both
$X_i$ and $Y_i$.
\end{itemize}
In particular, if the distribution $F$ of $X_1$ belongs to the
maximal domain of attraction of the Fr\'{e}chet distribution
$\Phi_{\alpha}$, then though the distribution $F_Y$ of $Y_1$ belongs
to the maximal domain of attraction of the Gumbel distribution, the
scaling factor depends on $\alpha$. This cannot happen in the i.i.d.
situation and, intuitively, it means that in the subordinated case
{\it the long range dependent sequence $\{X_i,i\ge 1\}$ also
contributes information to the asymptotic behavior of extreme sums}.

Moreover, we may have two LRD sequences $\{X_i,i\ge 1\}$,
$\{Y_i,i\ge 1\}$, the first one as in (\ref{model}), the second one
defined by $Y_i=G(X_i')$ with a sequence $\{X_i',i\ge 1\}$ defined
as in (\ref{model}), with the same covariance, with the same
marginals, but completely different
behavior of extremal terms. \\

Of course, it would be desirable to obtain some information about
limiting behaviour not only of extreme sums, but for sample maxima
as well. It should be pointed out that our method is not
appropriate. This is still an open problem to derive limiting
behaviour of maxima in the model (\ref{model}). In a different
setting, the case of stationary stable processes generated by
conservative flow, the
problem is treated in \cite{Samorodnitsky}.\\

We will use the following convention concerning integrals. If
$-\infty<a<b<\infty$ and $h$, $g$ are left-continuous and
right-continuous functions, respectively, then
$$
\int_a^bgdh=\int_{[a,b)}gdh \qquad \mbox{\rm and} \qquad
\int_a^bhdg=\int_{(a,b]}hdg ,
$$
whenever these integrals make sense as Lebesgue-Stjeltjes integrals.
The integration by parts formula yields
$$
\int_a^bg dh +\int_a^b hdg =h(b)h(b)-f(a)g(a) .
$$
We shall write $g\in RV_{\alpha}$ ($g\in SV$) if $g$ is regularly
varying at infinity with index $\alpha$ (slowly varying at
infinity).

In what follows $C$ will denote a generic constant which may be
different at each of its appearances. Also, for any sequences $a_n$
and $b_n$, we write $a_n\sim b_n$ if $\lim_{n\to\infty}a_n/b_n=1$.
Further, let $\ell(n)$ be a slowly varying function, possibly
different at each place it appears. On the other hand, $L(\cdot)$,
$L_0(\cdot)$, $L_1(\cdot)$, $L_1^*(\cdot)$, etc., are slowly varying
functions, fixed form the time they appear. Moreover, $g^{(k)}$
denotes the $k$th order derivative of a function $g$ and $Z$ is a
standard normal random variable. For any stationary sequence
$\{V_i,i\ge 1\}$, we will denote by $V$ the random variable with the
same distribution as $V_1$.
\section{Statement of results}\label{sec-results}
Let $F_{\epsilon}$ be the marginal distribution function of the
centered i.i.d. sequence $\{\epsilon_i,i\ge 1\}$. Also, for a given
integer $p$, the derivatives $F^{(1)}_{\epsilon}, \ldots,
F^{(p+3)}_{\epsilon}$ of $F_{\epsilon}$ are assumed to be bounded
and integrable. Note that these properties are inherited by the
distribution function $F$ of $X_1$ as well (cf. \cite{HoHsing} or
\cite{Wu2003}). Furthermore, assume that $\Exp \epsilon_1^4<\infty$.
These conditions are needed to
establish the reduction principle for the empirical process and will be assumed throughout the paper. \\

To study sums of $k_n$ largest observations, we shall consider the
following forms of $F$. For the statements below concerning regular
variation and domain of attractions we refer to \cite{deHaan},
\cite[Chapter 3]{EmbrechtsKluppelbergMikosch} or \cite{Horvath1987}.

The first assumption is that the distribution $F$ satisfies the
following Von-Mises condition:
\begin{equation}\label{eq-Q-RV}
\lim_{x\to\infty}\frac{xf(x)}{1-F(x)}=\alpha>0.
\end{equation}
Using notation from \cite{EmbrechtsKluppelbergMikosch}, the
condition (\ref{eq-Q-RV}) will be referred as $X\in
MDA(\Phi_{\alpha})$, since (\ref{eq-Q-RV}) implies that $X$ belongs
to the maximal domain of attraction of the Fr\'{e}chet distribution
with index $\alpha$. Then
\begin{equation}\label{eq-Q-RV-0}
Q(1-y)=y^{-1/\alpha}L_1(y^{-1}), \shortmbox{ as} y\to 0,
\end{equation}
and the density-quantile function $fQ(y)=f(Q(y))$ satisfies
\begin{equation}\label{eq-Q-RV-1}
fQ(1-y)=y^{1+1/\alpha}L_2(y^{-1}), \shortmbox{ as} y\to 0,
\end{equation}
where $L_2(u)=\alpha (L_1(u))^{-1}$.\\

The second type of assumption is that $F$ belongs to the maximal
domain of attraction of the double exponential Gumbel distribution,
written as $X\in MDA(\Lambda)$. Then the corresponding Von-Mises
condition implies
\begin{equation}\label{eq-Q-SV}
\lim_{y\to 0}\frac{fQ(1-y)\int_{1-y}^1(1-u)/fQ(u)du}{y^2}=1.
\end{equation}
Thus, with
$L_3(y^{-1})=\left(y^{-1}\int_{1-y}^1(1-u)/fQ(u)du\right)^{-1}$ one
has
$$
fQ(1-y)=yL_3(y^{-1}),
$$
and $L_3$ is slowly varying at infinity.\\

To study the effect of subordination, we will consider the
corresponding assumptions on $F_Y$, referred to later as $Y\in
MDA(\Phi_{\alpha_0})$ and $Y\in MDA(\Lambda)$, respectively:
\begin{equation}\label{eq-fQ-Y-RV-1}
Q_Y(1-y)=y^{-1/\alpha_0}L_1^*(y^{-1})\shortmbox{and}
f_YQ_Y(1-y)=y^{1+1/\alpha_0}L_2^*(y^{-1}), \shortmbox{ as} y\to 0,
\end{equation}
with $L_2^*(u)=\alpha_0 (L_1^*(u))^{-1}$, and
$$
f_YQ_Y(1-y)=yL_3^*(y^{-1}),
$$
where $L_3^*$ is defined in the corresponding way as $L_3$.\\

Recall that $Q_n(y)=\inf\{x:F_n(x)\ge y\}=X_{k:n}$ if
$\frac{k-1}{n}<y\le \frac{k}{n}$. Let
$T_n(m,k)=\sum_{i=m+1}^{n-k}Y_{i:n}$  and
$$ \mu_n(m,k)=n\int_{m/n}^{1-k/n}Q_Y(y)dy.
$$
The main result of this paper is the following theorem.
\begin{thm}\label{thm-extreme-sums}
Let $G(x)=Q_Y(F(x))$. Let $k_n=n^{\xi}$, where $\xi\in (0,1)$ is
such that
$$
\xi>\left\{\begin{array}{ll}
\frac{\beta+1/\alpha}{1+1/\alpha-1/\alpha_0},&
\shortmbox{if} X\in MDA(\Phi_{\alpha}), \;Y\in MDA(\Phi_{\alpha_0}), \quad (*)\\
\frac{\beta+1/\alpha}{1+1/\alpha},&
\shortmbox{if}X\in MDA(\Phi_{\alpha}), \;Y\in MDA(\Lambda),\quad (**)\\
\frac{\beta}{1-1/\alpha_0},&
\shortmbox{if}X\in MDA(\Lambda), \;Y\in MDA(\Phi_{\alpha_0}),\quad (***)\\
\beta,& \shortmbox{if}X\in MDA(\Lambda),\; Y\in MDA(\Lambda),\quad
(****).
\end{array} \right.
$$
Assume that $\Exp Y<\infty$. Let $p$ be the smallest positive
integer such that $(p+1)(2\beta-1)>1$ and assume that for
$r=1,\ldots,p$,
\begin{equation}\label{eq-cond-fQ-fQ}
\int_{1/2}^1F^{(r)}(Q(y))dQ_Y(y)=\int_{1/2}^1\frac{F^{(r)}(Q(y))}{f_YQ_Y(y)}dy<\infty.
\end{equation} Let
$$
A_n=\left\{\begin{array}{ll}
\left(\frac{n}{k_n}\right)^{1+1/\alpha-1/\alpha_0}L_{21}\left(\frac{n}{k_n}\right),&
\mbox{\rm if } X\in MDA(\Phi_{\alpha}), Y\in MDA(\Phi_{\alpha_0}),\\
\left(\frac{n}{k_n}\right)^{1+1/\alpha}L_{22}\left(\frac{n}{k_n}\right),&
\mbox{\rm if }X\in MDA(\Phi_{\alpha}), Y\in MDA(\Lambda),\\
\left(\frac{n}{k_n}\right)^{1-1/\alpha_0}L_{23}\left(\frac{n}{k_n}\right),&
\mbox{\rm if }X\in MDA(\Lambda), Y\in MDA(\Phi_{\alpha_0}),\\
\left(\frac{n}{k_n}\right)L_{24}\left(\frac{n}{k_n}\right),&
\mbox{\rm if }X\in MDA(\Lambda), Y\in MDA(\Lambda).
\end{array} \right.
$$
where $L_{21}, L_{22}, L_{23}, L_{24}$ are slowly varying functions
to be specified later on. Then
$$
A_n\sigma_{n,1}^{-1}\left(\sum_{j=n-k_n+1}^nY_{j:n}-n\int_{1-k_n/n}^1Q_Y(y)dy\right)\convdistr
Z.
$$
\end{thm}
The corresponding cases concerning assumptions on $X$ and $Y$ will
be referred as Case 1, Case 2, Case 3 and Case 4.
\begin{cor}
Under the conditions of Theorem {\rm\ref{thm-extreme-sums}}, if either $X\in MDA(\Phi_{\alpha})$ or
$X\in MDA(\Lambda)$, then
$$
\left(\frac{n}{k_n}\right)\ell(n)\sigma_{n,1}^{-1}\left(\sum_{j=n-k_n+1}^nY_{j:n}-n\int_{1-k_n/n}^1Q_Y(y)dy\right)\convdistr
Z.
$$
\end{cor}

In the subordinated case we have chosen to work with $G=Q_YF$ to
illustrate phenomena rather then deal with technicalities. One could
work with general functions $G$, but then one would need to assume
that $G$ has the power rank 1 (see \cite{HoHsing} for the
definition). Otherwise the scaling $\sigma_{n,1}^{-1}$ is not
correct. To see that $G(\cdot)=Q_YF(\cdot)$ has the power rank 1,
note that for $G_{\infty}(x):=\int_{-\infty}^{\infty}G(x+t)dF(t)$ we
have
$$
\frac{d}{dx}G_{\infty}(x)=\int_{-\infty}^{\infty}\frac{f(x+t)}{f_YQ_YF(x+t)}dF(t).
$$
Substituting $x=0$ and changing variables $y=F(t)$ we obtain
$$
\frac{d}{dx}G_{\infty}(x)|_{x=0}=\int_{0}^{1}\frac{fQ(y)}{f_YQ_Y(y)}dy\not=
0.
$$
Furthermore, we must assume that the distribution of $Y=G(X)$
belongs to the appropriate domain of attraction. For example, if
$X\in MDA(\Phi_{\alpha})$ and $Y_i=X_i^{\rho}$, $\rho>0$, then $Y\in
MDA(\Phi_{\alpha/\rho})$, provided that the map $x\to x^{\rho}$ is
increasing on ${\matR}$. Otherwise, if for example $\rho=2$, one
needs to impose conditions not only on the right tail of $X$, but on
the left one as well.

Nevertheless, to illustrate flexibility for the choice of $G$, let
$G(x)=\log (x^+)^{\alpha}$, $\alpha>0$. If $X\in
MDA(\Phi_{\alpha})$, then $Y=G(X)$ belongs to $MDA(\Lambda)$.
Further, since $\Exp X=0$, the quantile function $Q(u)$ of $X$ must
be positive for $u>u_0$ with some $u_0\in (0,1)$. Since the map
$x\to \log (x^+)^{\alpha}$ is increasing, $Q_Y(u)=Q_{\alpha\log
(X^+)}(u)=\alpha \log Q(u)$ for $u>u_0$. Consequently, from Theorem
\ref{thm-extreme-sums} we obtain the following corollary.
\begin{cor}\label{cor-logarithmic}
If (**) holds and $X\in MDA(\Phi_{\alpha})$, then
$$
\left(\frac{n}{k_n}\right)^{1+1/\alpha}L_{22}\left(\frac{n}{k_n}\right)\sigma_{n,1}^{-1}\left(\sum_{j=n-k_n+1}^n\log
(X_{j:n}^+)^{\alpha}-n\int_{1-k_n/n}^1\log Q(y)dy\right)\convdistr
Z.
$$

\end{cor}
\subsection{Remarks}
\begin{remark}\label{rem-infinitevariance}{\rm
From the beginning we assumed that $\Exp \epsilon_1^4<\infty$, thus,
in Cases 1 and 2 we have the requirement $\alpha\ge 4$ and this is
the only constrain on this parameter. Condition $\Exp Y<\infty$
requires $\alpha_0>1$ in case of $Y\in MDA (\Phi_{\alpha_0})$. In
view of (*), to be able to choose $\xi<1$ we need to have
$\alpha_0>(1-\beta)^{-1}>2$. The same restriction comes in Case 3.}
\end{remark}
\begin{remark}{\rm The conditions
(*)-(****) on $\xi$ are somehow restrictive. They come form the
quality of the rates in the reduction principle for the empirical
processes.
 }
\end{remark}
\begin{remark}{\rm
Appropriate results concerning the law of the iterated logarithm for
the extreme sums can be also stated, at least in the case of $Y\in
MDA(\Phi_{\alpha_0})$, by replacing $\sigma_{n,1}^{-1}$ in Theorems
\ref{thm-extreme-sums} with $\sigma_{n,1}^{-1}(\log\log n)^{-1/2}$.
In view of \cite{HaeuslerMason}, the most interesting phenomena
occur if $k_n$ is small ($k_n=o(\log\log n)$ in the i.i.d. case).
This, in view of the previous remark, cannot be treated in our
situation at all. }
\end{remark}
\begin{remark}{\rm
The conditions $D_r:=\int_{1/2}^1 F^{(r)}(Q(y))/f_YQ_Y(y)dy<\infty$
are not restrictive at all, since they are fulfilled for most
distributions with a regularly varying density-quantile function
$fQ(1-y)$, for those we refer to \cite{Parzen1979}. Consider for
example Case 1, and assume that the density $f$ is non-increasing on
some interval $[x_0,\infty)$. Then $F^{(r)}$ is regularly varying at
infinity with index $r+\alpha$. Thus, for some $x_1>x_0$
$$\int_{1/2\vee x_1}^1F^{(r)}(Q(y))/f_YQ_y(y)dy=\int_{1/2\vee x_1}^1(1-y)^{r/\alpha-1/\alpha_0}\ell(y)dy<\infty$$
for all $r\ge 1$ provided $\alpha_0>1$. If, additionally, we impose
the following {\it Cs\"{o}rg\H{o}-R\'{e}v\'{e}sz-type conditions}
(cf. \cite[Theorem 3.2.1]{CsorgoLN}):
\begin{itemize}
\item[{\rm (CsR1)}] $f$ exists on $(a,b)$, where
$a=\sup\{x:F(x)=0\}$, $b=\inf\{x:F(x)=1\}$, $-\infty\le a<b\le
\infty$,
\item[{\rm (CsR2)}] $\inf_{x\in (a,b)}f(x)>0$,
\end{itemize}
then in view of (CsR2) and the assumed boundness of derivatives
$F^{(r)}(\cdot)$, the integral $D_r$ is finite.
 }
\end{remark}
\begin{remark}\label{rem-1}{\rm
In the proof of Theorem \ref{thm-extreme-sums} we have to work with
both $Q(\cdot)$ and $fQ(\cdot)$. Therefore, we assumed the Von-Mises
condition (\ref{eq-Q-RV}) since it implies both (\ref{eq-Q-RV-0})
and (\ref{eq-Q-RV-1}). If one assumes only (\ref{eq-Q-RV-0}), then
(\ref{eq-Q-RV}) and, consequently, (\ref{eq-Q-RV-1}) hold, provided
a monotonicity of $f$ is assumed. Moreover, the von-Mises condition
is natural, since the existence of the density $f$ is explicitly
assumed.
 }
\end{remark}
\section{Proofs}
\subsection{Consequences of the reduction principle}
Let $p$ be a positive integer and let  \noindent
\begin{eqnarray*}
\lefteqn{S_{n,p}(x)=\sum_{i=1}^n(1_{\{X_i\le
x\}}-F(x))+\sum_{r=1}^p(-1)^{r-1}F^{(r)}(x)Y_{n,r}}\\
&=:&\sum_{i=1}^n(1_{\{X_i\le
x\}}-F(x))+V_{n,p}(x),\qquad\qquad\qquad\qquad
\end{eqnarray*}
where $F^{(r)}$ is the $r$th order derivative of $F$. Setting
$U_i=F(X_i)$ and $x=Q(y)$ in the definition of $S_{n}(\cdot)$, we
arrive at its uniform version,
\begin{eqnarray*}
\lefteqn{\tilde S_{n,p}(y)=\sum_{i=1}^n(1_{\{U_i\le
y\}}-y)+\sum_{r=1}^p(-1)^{r-1}F^{(r)}(Q(y))Y_{n,r}}\\
&=:&\sum_{i=1}^n(1_{\{U_i\le y\}}-y)+\tilde
V_{n,p}(y).\qquad\qquad\qquad\qquad
\end{eqnarray*}
Denote
$$
d_{n,p}=\left\{\begin{array}{ll} n^{-(1-\beta)}L_0^{-1}(n)(\log n)^{5/2}(\log\log n)^{3/4}, & \frac{p+1}{2\beta-1}\ge 1\\
n^{-p(\beta-\frac{1}{2})}L_0^p(n)(\log n)^{1/2}(\log\log n)^{3/4}, &
\frac{p+1}{2\beta-1}<1
\end{array}\right. .
$$

We shall need the following lemma, referred to as the reduction
principle.
\begin{lem}[\cite{Wu2003}]\label{thm-HoHsing}
Let $p$ be a positive integer. Then, as $n\to\infty$,
$$
\Exp \sup_{x\in {\matR}}\left|\sum_{i=1}^n(1_{\{X_i\le
x\}}-F(x))+\sum_{r=1}^p(-1)^{r-1}F^{(r)}(x)Y_{n,r}\right|^2=O(\Xi_n+n(\log
n)^2),
$$
where
$$
\Xi_n=\left\{\begin{array}{ll} O(n), & (p+1)(2\beta-1)>1\\
O(n^{2-(p+1)(2\beta-1)}L_0^{2(p+1)}(n)), & (p+1)(2\beta-1)<1
\end{array}\right. .
$$
\end{lem}

Using Lemma \ref{thm-HoHsing} we obtain (cf.
\cite{CsorgoKulik2006a})
\begin{eqnarray*}
\lefteqn{\sigma_{n,p}^{-1}\sup_{x\in
{\matR}}|S_{n}(x)|}\\
&=&\left\{\begin{array}{ll} O_{a.s}(n^{-(\frac{1}{2}-p(\beta-\frac{1}{2}))}L_0^{-p}(n)(\log n)^{5/2}(\log\log n)^{3/4}), &\frac{p+1}{2\beta-1}>1\\
O_{a.s}(n^{-(\beta-\frac{1}{2})}L_0(n)(\log n)^{1/2}(\log \log
n)^{3/4}), & \frac{p+1}{2\beta-1}<1
\end{array}\right. .
\end{eqnarray*}
Since (see (\ref{eq-varance-behaviour}))
$$
\frac{\sigma_{n,p}}{\sigma_{n,1}}\sim
n^{-(\beta-\frac{1}{2})(p-1)}L_0^{p-1}(n)
$$
we obtain
\begin{eqnarray*}
\lefteqn{\sup_{x\in {\matR}}|\beta_n(x)+\sigma_{n,1}^{-1}V_{n,p}(x)|=}\\
& = & \frac{\sigma_{n,p}}{\sigma_{n,1}}\sup_{x\in
{\matR}}\left|\sigma_{n,p}^{-1}\sum_{i=1}^n(1_{\{X_i\le
x\}}-F(x))+\sigma_{n,p}^{-1}V_{n,p}(x)\right|=o_{a.s}(d_{n,p}).
\end{eqnarray*}
Consequently, via $\{\alpha_n(y),y\in (0,1)\}=\{\beta_n(Q(y)),y\in
(0,1)\}$,
\begin{equation}\label{approx-unif-empirical}
\sup_{y\in (0,1)}|\alpha_n(y)+\sigma_{n,1}^{-1}\tilde
V_{n,p}(y)|=O_{a.s}(d_{n,p}).
\end{equation}
We have
\begin{eqnarray}
\lefteqn{ A_n\sigma_{n,1}^{-1}\int_{1-a_n/n}^{1-1/n}\tilde
V_{n,p}(y)dQ_Y(y)=A_n\sigma_{n,1}^{-1}\int_{1-a_n/n}^{1-1/n}\frac{\tilde
V_{n,p}(y)}{f_YQ_Y(y)}dy}\label{eq-technical-4}\\
&=&-\left(A_n\int_{1-a_n/n}^{1-1/n}\frac{fQ(y)}{f_YQ_Y(y)}dy\right)\left[\left(\sigma_{n,1}^{-1}\sum_{i=1}^nX_i\right)+
o_P(\sigma_{n,1}^{-1})\right].\nonumber
\end{eqnarray}
Let $$L_{11}(u)=L^*_2(u)/L_2(u),\qquad
L_{21}(u)=(1/\alpha-1/\alpha_0+1)L_{11}(u),$$
$$L_{12}(u)=L^*_3(u)/L_2(u),\qquad
L_{22}(u)=(1/\alpha+1)L_{12}(u),$$
$$L_{13}(u)=L^*_2(u)/L_3(u),\qquad
L_{23}(u)=(-1/\alpha+1)L_{13}(u),$$
$$L_{14}(u)=L^*_3(u)/L_3(u),\qquad L_{24}(u)=L_{14}(u).$$
\begin{lem}\label{lem-convergence-full-sums}
Let $p$ be a positive integer. Assume that for $r=1,\ldots,p$,
{\rm(\ref{eq-cond-fQ-fQ})} holds. Then
$$
A_n\sigma_{n,1}^{-1}\int_{1-k_n/n}^{1-1/n}\tilde
V_{n,p}(y)dQ_Y(y)\convdistr Z.
$$
\end{lem}
{\it Proof.} In view of (\ref{eq-technical-4}), we need only to
study the asymptotic behavior, as $n\to\infty$, of
$A_n\int_{1-k_n/n}^{1-1/n}\frac{fQ(y)}{f_YQ_Y(y)}dy=:A_nK_n$ and to show that $A_nK_n\sim 1$.\\

We have by Karamata's Theorem:

\noindent In Case 1,
\begin{eqnarray*}
\lefteqn{K_n=\int_{1-k_n/n}^{1-1/n}(1-y)^{1/\alpha-1/\alpha_0}\left(L_{11}((1-y)^{-1})\right)^{-1}dy}\\
&\sim&
(1/\alpha-1/\alpha_0+1)^{-1}\left(\frac{k_n}{n}\right)^{1+1/\alpha-1/\alpha_0}\left(L_{11}\left(\frac{n}{k_n}\right)\right)^{-1}\\
&\sim&
\left(\frac{k_n}{n}\right)^{1+1/\alpha-1/\alpha_0}\left(L_{21}\left(\frac{n}{k_n}\right)\right)^{-1}.
\end{eqnarray*}
\noindent In Case 2,
\begin{eqnarray*}
\lefteqn{K_n=\int_{1-k_n/n}^{1-1/n}(1-y)^{1/\alpha}\left(L_{12}((1-y)^{-1})\right)^{-1}dy}\\
&\sim&
(1/\alpha+1)^{-1}\left(\frac{k_n}{n}\right)^{1+1/\alpha}\left(L_{12}\left(\frac{n}{k_n}\right)\right)^{-1}\sim
\left(\frac{k_n}{n}\right)^{1+1/\alpha}\left(L_{22}\left(\frac{n}{k_n}\right)\right)^{-1}.
\end{eqnarray*}
\noindent In Case 3,
\begin{eqnarray*}
\lefteqn{K_n=\int_{1-k_n/n}^{1-1/n}(1-y)^{-1/\alpha}\left(L_{13}((1-y)^{-1})\right)^{-1}dy}\\
&\sim&
(-1/\alpha+1)^{-1}\left(\frac{k_n}{n}\right)^{1-1/\alpha}\left(L_{13}\left(\frac{n}{k_n}\right)\right)^{-1}\sim
\left(\frac{k_n}{n}\right)^{1-1/\alpha}\left(L_{23}\left(\frac{n}{k_n}\right)\right)^{-1}.
\end{eqnarray*}
\noindent In Case 4,
\begin{eqnarray*}
\lefteqn{K_n=\int_{1-k_n/n}^{1-1/n}\left(L_{14}((1-y)^{-1})\right)^{-1}dy}\\
&\sim&
(-1/\alpha+1)^{-1}\left(\frac{k_n}{n}\right)\left(L_{14}\left(\frac{n}{k_n}\right)\right)^{-1}\sim
\left(\frac{k_n}{n}\right)\left(L_{14}\left(\frac{n}{k_n}\right)\right)^{-1}.
\end{eqnarray*}
Thus, in either case, $A_nK_n\sim 1$. \koniec
\begin{lem}\label{lem-1}
For any $k_n\to\infty$, $k_n=o(n)$
$$
\frac{U_{n-k_n:n}}{1-k_n/n}\convprob 1.
$$
\end{lem}
{\it Proof.} In view of (\ref{approx-unif-empirical}) one obtains
$$
\sup_{y\in (0,1)}|u_n(y)|=\sup_{y\in (0,1)}|\alpha_n(y)|=O_P(1).
$$
Consequently,
\begin{eqnarray*}
\sup_{y\in (0,1)}|y-U_n(y)|&=&\sup_{y\in
(0,1)}\sigma_{n,1}n^{-1}|u_n(y)|=\sup_{y\in
(0,1)}\sigma_{n,1}n^{-1}|\alpha_n(y)|\nonumber\\
&=&O_{P}(\sigma_{n,1}n^{-1}).
\end{eqnarray*}
Thus, the result follows by noting that
$U_n(1-k_n/n)=U_{n-k_n:n}$.\koniec

An easy consequence of (\ref{approx-unif-empirical}) is the
following result.
\begin{lem}\label{lem-2}
For any $k_n\to 0$,
$$
\sup_{y\in
(1-k_n/n,1)}|\alpha_n(y)|=O_{a.s.}(d_{n,p})+O_{P}(f(Q(1-k_n/n))).
$$
\end{lem}
\subsection{Proof of Theorem \ref{thm-extreme-sums}}
To obtain the limiting behavior of sums of extremes, we shall use
the following decomposition: Since $E_n(\cdot)$ has no jumps after
$U_{n:n}$ and $Y_j=Q_YF(X_j)=Q_Y(U_j)$, we have
\begin{eqnarray*}
\lefteqn{A_n\sigma_{n,1}^{-1}\left(\sum_{j=n-k_n+1}^nY_{j:n}-n\int_{1-k_n/n}^1Q_Y(y)dy\right)}\nonumber\\
&= &
A_n\sigma_{n,1}^{-1}\left(\sum_{j=n-k_n+1}^nQ_Y(U_{j:n})-n\int_{1-k_n/n}^1Q_Y(y)dy\right)\nonumber
\\
&= & A_n\sigma_{n,1}^{-1}\left(n\int
_{U_{n-k_n:n}}^{U_{n:n}}Q_Y(y)dE_n(y)-n\int_{1-k_n/n}^1Q_Y(y)dy\right)\nonumber\\
&= & A_n\sigma_{n,1}^{-1}\left(n\int
_{U_{n-k_n:n}}^{1}Q_Y(y)dE_n(y)-n\int_{1-k_n/n}^1Q_Y(y)dy\right)\nonumber\\
& = & A_n\sigma_{n,1}^{-1}n\left\{\truncatedintegral
(y-E_n(y))dQ_Y(y)\right.\nonumber\\
&&\left.+\int_{1-\frac{1}{n}}^1(y-E_n(y))dQ_Y(y)+\int_{U_{n-k_n:n}}^{1-k_n/n}(1-\frac{k_n}{n}-E_n(y))dQ_Y(y)\right\}\nonumber\\
&= & -A_n\truncatedintegral
\alpha_n(y)dQ_Y(y)-A_n\int_{1-1/n}^1\alpha_n(y)dQ_Y(y)\nonumber
\\&&+A_n\sigma_{n,1}^{-1}n\int_{U_{n-k_n:n}}^{1-k_n/n}(1-\frac{k_n}{n}-E_n(y))dQ_Y(y)=:I_1+I_2+I_3.\label{eq-sum-decomposition}
\end{eqnarray*}
We will show that $I_1$  yields the asymptotic normality.
Further, we will show that the latter two integrals are
asymptotically negligible.

Each term will be treated in a separate section. Let $p$ be the
smallest integer such that $(p+1)(2\beta-1)>1$, so that
$d_{n,p}=n^{-(1-\beta)}\ell(n)$.
\subsubsection{First term}
Let $\psi_{\mu}(y)=(y(1-y))^{\mu}$, $y\in [0,1]$, $\mu>0$.\\

For $k_n=n^{\xi}$ and arbitrary small $\delta>0$ one has by
(\ref{approx-unif-empirical}),
\begin{eqnarray*}
\lefteqn{A_n\left|\alpha_n(y)+\sigma_{n,1}^{-1}\tilde
V_{n,p}(y)\right|=O_{a.s}\left(A_nd_{n,p}\right)}\label{eq-technical-2}\\
&=& \left\{\begin{array}{ll}
n^{-(\xi+\xi/\alpha-\xi/\alpha_0-1/\alpha+1/\alpha_0-\beta-\delta)}
,&
\shortmbox{if} X\in MDA(\Phi_{\alpha}), Y\in MDA(\Phi_{\alpha_0}),\\
n^{-(\xi+\xi/\alpha-1/\alpha-\beta-\delta)} ,&
\shortmbox{if}X\in MDA(\Phi_{\alpha}), Y\in MDA(\Lambda),\\
n^{-(\xi-\xi/\alpha_0+1/\alpha_0-\beta-\delta)},&
\shortmbox{if}X\in MDA(\Lambda), Y\in MDA(\Phi_{\alpha_0}),\\
n^{-(\xi-\beta-\delta)},& \shortmbox{if}X\in MDA(\Lambda), Y\in
MDA(\Lambda).
\end{array} \right.
\nonumber
\end{eqnarray*}
Let
$$
J_n=A_n\left|\truncatedintegral
\frac{\left|\alpha_n(y)+\sigma_{n,1}^{-1}\tilde
V_{n,p}(y)\right|}{\psi_{\mu}(y)}\psi_{\mu}(y)dQ_Y(y)\right| .
$$
Case 1: Since condition (*) on $\xi$ holds,
$$
1/\alpha_0<\xi+\xi(1/\alpha-1/\alpha_0)-1/\alpha+1/\alpha_0-\beta .
$$
Set $\mu=(\alpha_0-\delta)^{-1}$ with $\delta>0$ so small that
$$
\mu<\xi+\xi(1/\alpha-1/\alpha_0)-1/\alpha+1/\alpha_0-\beta-\delta .
$$
Then, we have $\Exp (Y^+)^{1/\mu+\delta/2}<\infty$. The latter
condition is sufficient for the finiteness of
$\int_{x_1}^1\psi_{\mu}(y)dQ_Y(y)$, where $x_1=\inf\{y:Q_Y(y)\ge
0\}$, (see \cite[Remark
2.4]{ShaoYu}). Thus,
$$
J_n=o_{a.s}(A_nd_{n,p}n^{\mu})\int_{x_1}^1\psi_{\mu}(y)dQ_Y(y)=o_{a.s}(1)O(1).
$$
Since in Case 3, (***) holds, a similar approach yields that in
this case $J_n=o_{a.s}(1)$.\\

\noindent Case 2: If $Y\in MDA(\Lambda)$ then $\Exp
(Y^+)^{\alpha_0}<\infty$ for all $\alpha>0$ (see \cite[Corollary
3.3.32]{EmbrechtsKluppelbergMikosch}). Thus, in view of (**), choose
arbitrary small $\delta>0$ and $\alpha_0$ so big that $\Exp
(Y^+)^{\alpha_0}<\infty$ and
$$
\frac{1}{\alpha_0-\delta}<\xi+\xi/\alpha-1/\alpha-\beta-\delta .
$$
Set $\mu=(\alpha_0-\delta)^{-1}$ and continue as in the Case 1. A
similar reasoning applies to Case 4, provided $\xi>\beta$. Thus, in
either case
$$
A_n\left|\truncatedintegral
\left(\alpha_n(y)+\sigma_{n,1}^{-1}\tilde
V_{n,p}(y)\right)dQ_Y(y)\right|=o_{a.s}(1) \label{eq-technical-3}.
$$
Now, the asymptotic normality of $I_1$ follows from Lemma \ref{lem-convergence-full-sums}.
\subsubsection{Second term}
We have
\begin{eqnarray*}
\lefteqn{A_n\int_{1-1/n}^1\alpha_n(y)dQ_Y(y)}\\
&= &
-A_n\sigma_{n,1}^{-1}n\int_{1-1/n}^1(1-E_n(y))dQ_Y(y)+A_n\sigma_{n,1}^{-1}n\int_{1-1/n}^1(1-y)dQ_Y(y)\\
&:=&J_1+J_2.
\end{eqnarray*}
Since $\Exp J_1=J_2$, it suffices to show that $J_2=o(1)$.\\

\noindent Case 1: We have by Karamata's Theorem
\begin{eqnarray*}
\lefteqn{J_2=A_n\sigma_{n,1}^{-1}n\int_{1-1/n}^1\frac{(1-y)}{(1-y)^{1+1/\alpha_0}L_2^*(y^{-1})}dy}\\
&\sim&
\left(\frac{n}{k_n}\right)^{1+1/\alpha-1/\alpha_0}n^{\beta-3/2}n(\frac{1}{n})^{1-1/\alpha_0}\ell(n)\ell(n/k_n)
\end{eqnarray*}
which converges to 0 using the assumption (*).\\

\noindent Likewise, in Case 3,
$$
J_2\sim
\left(\frac{n}{k_n}\right)^{1-1/\alpha_0}n^{\beta-3/2}n(\frac{1}{n})^{1-1/\alpha_0}\ell(n)\ell(n/k_n)
$$
which converges to 0 using the assumptions (***).\\

\noindent Case 2: We have,
\begin{eqnarray*}
\lefteqn{J_2=A_n\sigma_{n,1}^{-1}n\int_{1-1/n}^1\frac{1-y}{f_YQ_Y(y)}dy}\\
&\sim& A_n\sigma_{n,1}^{-1}n^{-1}(f_YQ_Y(1-1/n))^{-1}\sim
A_n\sigma_{n,1}^{-1}\ell(n)\ell(n/k_n)
\end{eqnarray*}
which converges to 0, using the assumption (**). The same argument
applies to Case 4. Therefore, in either case, $I_2=o_P(1)$.
\subsubsection{Third term}
To prove that $I_3=o_P(1)$, let $y$ be in the interval with the
endpoints $U_{n-k_n:n}$ and $1-k_n/n$. Then
$$
\left|1-E_n(y)-\frac{k_n}{n}\right|\le |E_n(1-k_n/n)-(1-k_n/n)|.
$$
\noindent Case 1: By Lemma \ref{lem-1} and $Y\in
MDA(\Phi_{\alpha_0})$, we have
\begin{equation}\label{eq-Q-Q-to-1}
Q_Y(1-k_n/n)/Q_Y(U_{n-k_n:n})\convprob 1.\end{equation} Hence, by
condition (*),
\begin{eqnarray}\label{eq-technical-1}
\lefteqn{\left(\frac{n}{k_n}\right)^{1+1/\alpha-1/\alpha_0}\ell(n/k_n)Q_Y(1-k_n/n)d_{n,p}}\nonumber\\
&=&n^{1+1/\alpha}\ell(n)\ell(n/k_n)n^{-\xi(1+1/\alpha)}d_{n,p}\to 0.
\end{eqnarray}
Also, by (\ref{eq-Q-RV-1}) and (\ref{eq-fQ-Y-RV-1}),
\begin{equation}\label{eq-technical-1a}
A_nQ_Y(1-k_n/n)fQ(1-k_n/n)\sim C
L_{21}\left(\frac{n}{k_n}\right)\frac{L_1^*(n/k_n)}{L_1(n/k_n)}\sim
C
\end{equation}
Thus, by (\ref{eq-Q-Q-to-1}), (\ref{eq-technical-1}),
(\ref{eq-technical-1a}) and Lemma \ref{lem-2}
\begin{eqnarray*}
\lefteqn{I_3\le A_nQ_Y(1-k_n/n)|\alpha_n(1-k_n/n)|\frac{|Q_Y(1-k_n/n)-Q_Y(U_{n-k_n:n})|}{Q_Y(1-k_n/n)}}\\
&= & A_nQ_Y(1-k_n/n)\alpha_n(1-k_n/n)o_{p}(1)\\
&=&o_p\left(A_nQ_Y(1-k_n/n)fQ(1-k_n/n)\right)+o_p\left(A_nQ(1-k_n/n)d_{n,p}\right)=o_P(1).
\end{eqnarray*}
\noindent Case 3: By Lemma \ref{lem-1} and $Y\in
MDA(\Phi_{\alpha_0})$ we have (\ref{eq-Q-Q-to-1}). Since
$\xi>\beta$,
\begin{equation}\label{eq-technical-1-1}
\left(\frac{n}{k_n}\right)^{1-1/\alpha_0}\ell(n/k_n)Q_Y(1-k_n/n)d_{n,p}=n^{\beta-\xi}\to
0.
\end{equation}
Also, by (\ref{eq-Q-SV}) and (\ref{eq-fQ-Y-RV-1}),
\begin{eqnarray}\label{eq-technical-1a-1}
\lefteqn{\left(\frac{n}{k_n}\right)^{1-1/\alpha_0}L_{23}\left(\frac{n}{k_n}\right)Q_Y(1-k_n/n)fQ(1-k_n/n)}\hspace*{3cm}\nonumber\\
&\sim &C
L_{23}\left(\frac{n}{k_n}\right)\frac{L_3(n/k_n)}{L_2^*(n/k_n)}\sim
C.
\end{eqnarray}
Thus, by (\ref{eq-technical-1-1}), (\ref{eq-technical-1a-1}), we
conclude as above that $I_3=o_P(1)$.\\

\noindent Cases 2 and 4:
$$
T_n(\lambda)=A_n|\alpha_n(1-k_n/n)|\left|Q_Y(r_n^{+}(\lambda))-Q_Y(r_n^{-}(\lambda))\right|,
$$
where $r_n^{+}(\lambda)=1-\frac{k_n}{\lambda n}$,
$r_n^{-}(\lambda)=1-\frac{k_n}{\lambda n}$ and $1<\lambda<\infty$ is
arbitrary. Applying an argument as in the proof of Theorem 1 in
\cite{CsorgoHorvathMason}, we have
$$
\liminf_{n\to\infty}P(|I_3|<|T_n(\lambda)|)\ge
\lim\inf_{n\to\infty}P(r_n^{-}(\lambda)\le U_{n-k_n:n}\le
r_n^{+}(\lambda)).
$$
In view of Lemma \ref{lem-1}, the lower bound is 1. Thus,
$\lim_{n\to\infty}P(|I_3|<|T_n(\lambda)|)=1$. Further, by Lemma 4 in
\cite{Lo},
$$
\lim_{n\to\infty}(Q_Y(r_n^{+}(\lambda))-Q_Y(r_n^{-}(\lambda))){L_3^*(n/k_n)}=-\log
\lambda .
$$
Thus, for large $n$,
\begin{eqnarray*}
T_n(\lambda)& = &
A_n|\alpha_n(1-k_n/n)|(L_3^*(n/k_n))^{-1}|Q_Y(r_n^{+}(\lambda))-Q_Y(r_n^{-}(\lambda))|L_3^*(n/k_n)\\
&\le & C_1\frac{A_n}{L_3^*(n/k_n)}fQ(1-k_n/n)(\log
\lambda)+C_2\frac{A_n}{L_3^*(n/k_n)}d_{n,p}\log \lambda
\end{eqnarray*}
almost surely with some constants $C_1,C_2$. The second term, for
arbitrary $\lambda$, converges to 0 by the choice of $\xi$. Also,
$$
A_n\frac{fQ(1-k_n/n)}{L_3^*(n/k_n)}\le\left\{\begin{array}{ll}
\left(\frac{n}{k_n}\right)^{1+1/\alpha}L_{22}\left(\frac{n}{k_n}\right)\frac{\left(\frac{k_n}{n}\right)^{1+1/\alpha}L_2\left(\frac{n}{k_n}\right)}{L_3^*\left(\frac{n}{k_n}\right)},&
\mbox{\rm in Case 2},\\
\left(\frac{n}{k_n}\right)L_{24}\left(\frac{n}{k_n}\right)\frac{\left(\frac{k_n}{n}\right)L_3\left(\frac{n}{k_n}\right)}{L_3^*\left(\frac{n}{k_n}\right)},&
\mbox{\rm in Case 4}.
\end{array}
\right.
$$
In either case, the above expressions are asymptotically equal to 1.
Thus, we have for sufficiently large $n$, $T_n(\lambda)\le
C_1\log\lambda$ almost surely. Thus,
$\lim_{n\to\infty}P(|T_n(\lambda)|\le C_1\log\lambda)=1$.
Consequently,
\begin{eqnarray*}
\lefteqn{\lim_{n\to\infty}P(|I_3|> C_1\log\lambda)=}\\
& = & \lim_{n\to\infty}P(|I_3|> C_1\log\lambda,|T_n(\lambda)|\le C_1
\log\lambda)+\lim_{n\to\infty}P(|T_n(\lambda)|> C_1\log\lambda)\\
&\le & \lim_{n\to\infty}P(|I_3|> |T_n(\lambda)|)+0=0
\end{eqnarray*}
and thus $I_3=o_P(1)$.
 \koniec
\ack This work was done during my stay at Carleton University. I am
thankful to Professors Barbara Szyszkowicz and Miklos Cs\"{o}rg\H{o}
for the support and helpful remarks.

\end{document}